\magnification=1200

\centerline
{\bf Hasse Principle for $G$--quadratic forms}

\bigskip
\centerline 
{Eva Bayer--Fluckiger, Nivedita Bhaskhar and Raman Parimala}

\bigskip
\bigskip
{\bf Introduction.}

\bigskip
Let $k$ be a global field of characteristic $\not = 2$. The classical Hasse--Minkowski
theorem states that if two quadratic forms become isomorphic over all the completions
of $k$, then they are isomorphic over $k$ as well. It is natural to ask whether this
is true for  {\it $G$--quadratic forms},  where $G$ is a finite group. In the case of
number fields the Hasse principle for $G$--quadratic forms does not hold in general, as shown by J. Morales [M 86]. The aim of the present
paper is to study this question when $k$ is a global field of positive characteristic.
We give a sufficient criterion for the Hasse principle to hold (see th.  2.1.), and also
give counter--examples. These counter--examples are of a different nature than
those for number fields : indeed, if $k$ is a global field of positive characteristic,
then the Hasse principle does hold for $G$--quadratic forms on projective $k[G]$--modules
(see cor. 2.3), and in particular if $k[G]$ is semi--simple, then the Hasse principle is
true for $G$--quadratic forms, contrarily to what happens in the case of number fields.
On the other hand, there are counter--examples in the non semi--simple case, as
shown in \S 3. Note that the Hasse principle holds in all generality
for $G$--trace forms (cf. [BPS 13]). 
\medskip
The third named author is partially supported by National Science Foundation grant DMS-1001872.
\bigskip
{\bf \S 1. Definitions, notation and basic facts}

\bigskip
Let $k$ be a field of characteristic $\not = 2$. All modules are supposed to be left
modules.

\medskip
{\it $G$--quadratic spaces}
\medskip
Let $G$ be a finite group, and let
$k[G]$ be the associated group ring. 
A {\it $G$--quadratic space} is a pair $(V,q)$, where $V$ is a $k[G]$--module
that is a finite dimensional $k$--vector space, 
and $q : V \times V \to k$ is a non--degenerate symmetric bilinear
form such that $$q(gx,gy) = q(x,y)$$ for all $x,y \in V$ and all $g \in G$. 
\medskip
Two $G$--quadratic spaces  $(V,q)$ and $(V',q')$ are {\it isomorphic} if there exists
an isomorphism of $k[G]$--modules $f : V \to V'$ such that $q'(f(x),f(y)) = q(x,y)$ for
all $x,y \in V$. If this is the case, we write $(V,q) \simeq_G (V',q')$, or simply $q \simeq_G q'.$

\bigskip

{\it Hermitian forms }

\medskip

Let $R$ be a ring endowed with an involution $r \mapsto \overline r$. For any 
$R$--module $M$, we denote by $M^*$ its dual  ${\rm Hom}_R(M,R)$. Then $M^*$  has an $R$--module structure given by $(rf)(x) =  f(x) \overline r$ for all $r \in R$, $x \in M$ and
$f \in M^*$. If $M$ and $N$ are two $R$--modules and if $f : M \to N$ is
a homomorphism of $R$--modules, then $f$ induces a homomorphism
$f^* : N^* \to M^*$ defined by $f^*(g) = g f$ for all $g \in N^*$, called the {\it adjoint} of $f$.

\medskip
A {\it hermitian form} is a pair $(M,h)$ where
$M$ is an $R$--module and $h : M \times M \to R$ is biadditive, satisfying the following two conditions:

\smallskip
 (1.1) $h(rx,sy) = r h(x,y) \overline s$ and $\overline {h(x,y)} = h(y,x)$ for
all $x,y \in M$ and all $r,s \in R$.

(1.2) The homomorphism $ h : M \to M^*$ given by  $y \mapsto h(\ , y )$ is an isomorphism.

\smallskip
Note that the existence of $h$ implies that $M$ is {\it self-dual}, i.e. isomorphic to its dual.

\medskip
  
 If $G$ is a finite group, then the group
  algebra  $R = k[G]$ has a natural $k$--linear involution, characterized by
  the formula $\overline g = g^{-1}$ for every $g \in G$. We have the following dictionary (see for
  instance [BPS 13, 2.1, Example]
  
   a) $R$--module $M  \ \Longleftrightarrow  k$-module $M$ with a $k$--linear action of $G$;
   
   b) $R$--dual $M^*\ \Longleftrightarrow k$--dual of $M$, with the contragredient (i.e. dual) action of $G$. 
   
  c) hermitian space $(M,h)$\ $\Longleftrightarrow$ 
  symmetric bilinear form on $M$, which is $G$--invariant and defines an isomorphism of $M$ onto its $k$--dual.
  
  \medskip
  Therefore a hermitian space over $k[G]$ corresponds to a $G$--quadratic space,
  as defined above. 

\bigskip

{\it Hermitian elements}

\medskip
 Let $E$ be a ring with an involution $\sigma : E \to E$ 
 and put 
$$E^{0} = \{ z \in E^{\times} \ | \ \sigma(z) =  z \}.$$  
If $z \in E^{0}$, the map $h_z : E \times E \to E$ defined by
$h_z(x,y) = x.z.\sigma(y)$ is a hermitian space over $E$;
conversely, every hermitian space over $E$ with underlying module $E$ is isomorphic to $h_z$
for some $z \in E^{0}$.
\medskip
Define an equivalence
relation on $E^{0} $ by setting 
$z \equiv z'$ if there exists $e \in E^\times$ with $z' = \sigma(e) z e$; this is equivalent to $ (E,h_z) \simeq (E,h_{z'})$. Let
$ H(E,\sigma)$
be the quotient of $E^{0}$ by this equivalence relation.
If $z \in E^{0} $, we denote by $[z]$ its class in $H(E,\sigma)$.

\bigskip

{\it Classifying hermitian spaces via hermitian elements}

\medskip
Let $(M,h_0)$ be 
a hermitian space over $R$. Set $E_M = {\rm End}(M)$. Let $\tau : E_M \to E_M$ be the involution of $E_M$ {\it induced} by $h_0$, i.e. 

$$\tau (e) = h_0^{-1} e^* h_0,  \quad {\rm for} \ \  e \in E_M, $$

\noindent where $e^*$ is the adjoint of  $e$. If $(M,h)$ is a hermitian space (with the same underlying module $M$), we have
$\tau(h_0^{-1} h) = h_0^{-1}(h_0^{-1}h)^*h_0 = h_0^{-1}h^*(h_0^{-1})^*h_0 = h_0^{-1} h$.
Hence $h_0^{-1} h$ is a hermitian element of $(E_M,\tau)$; let $[h_0^{-1} h]$ be its class in
$H(E_M,\tau).$

\bigskip
\noindent
{\bf Lemma 1.1.} (see for instance [BPS 13, lemma 3.8.1]) {\it Sending a hermitian space $(M,h)$ to the element $[h_0^{-1} h]$ of
$H(E_M,\tau)$ induces a bijection between the set of isomorphism classes
of hermitian spaces $(M,h)$ and the set $H(E_M,\tau)$. }

\bigskip
{\it Components of algebras with involution}
\medskip
Let $A$ be a finite dimensional $k$--algebra, and let 
$\iota : A \to A$ be a $k$--linear involution.
Let $R_A$ be the radical of $A$. Then $A/R_A$ is
a semi--simple $k$--algebra, hence we have a decomposition 
$A/R_A = \prod_{i = 1,\dots,r} M_{n_i}(D_i)$,
where $D_1,\dots,D_r$ are
division algebras. Let us denote by $K_i$ the center of $D_i$, and let $D^{\rm op}_i$ be
the opposite algebra of $D_i$. 
\medskip
Note that $\iota(R_A) = R_A$, hence $\iota$ induces an involution $\iota : A/R_A  \to A/R_A$. 
Therefore $A/R_A$ decomposes into a product of involution invariant factors. These can
be of two types~: either an involution invariant matrix algebra  $M_{n_i}(D_i)$, or a product
$M_{n_i}(D_i )\times M_{n_i}(D^{op}_i)$, with $M_{n_i}(D_i)$ and $M_{n_i}(D^{op}_i)$ exchanged
by the involution.  We say that a factor is {\it unitary} if the restriction of the involution to its
center is not the identity : in other words, either an involution invariant $M_{n_i}(D_i)$ with
$\iota | K_i$ not the identity, or a product $M_{n_i}(D_i) \times M_{n_i}(D^{op}_i)$ . Otherwise,
the factor is said to be of the first kind. In this case, the component is of the form 
$M_{n_i}(D_i)$ and the restriction of $\iota$ to $K_i$ is the identity. We say that the
component is {\it orthogonal} if after base change to a separable closure $\iota$ is given
by the transposition, and {\it symplectic} otherwise.
A component $M_{n_i}(D_i)$  
 is said to be  {\it split} if $D_i$ is a commutative field. 
 
 \bigskip
 {\it Completions}
 \medskip
 If $k$ is a global field and if $v$ is a place of $k$, we denote by $k_v$ the completion
 of $k$ at $v$. For any $k$--algebra $E$, set $E_v = E \otimes_k k_v$.  
 If $K/k$ is a field extension of finite degree and if $w$ is a place of $K$ above $v$,
 then we use the notation $w | v$. 

\bigskip
{\bf \S 2. Hasse principle}
\bigskip In this section, $k$ will be a global field of characteristic $\not = 2$. Let us
denote by $\Sigma_k$ the set of all places of $k$. 
The aim of this section is to give a sufficient criterion for the Hasse principle for
$G$--quadratic forms to hold. 
All modules are left modules, and finite dimensional $k$--vector spaces.

\bigskip
 \noindent
 {\bf Theorem 2.1.} {\it  Let $V$ be a $k[G]$--module ,
  and let $E = {\rm End}(V)$.  Let $R_E$ be the radical of $E$, and set
 $\overline E = E / R_E$. Suppose that all the orthogonal components of $\overline E$ are split, and let
 $(V,q)$, $(V,q')$  be two $G$--forms. Then $q \simeq_G q'$ over $k$ if and only if $q \simeq_G q'$
 over all the completions of $k$.}
 
\bigskip This is announced in [BP 13], and replaces th. 3.5 of [BP 11]. 
The proof of th. 2.1 relies on the following proposition
\bigskip
\noindent
{\bf Proposition 2.2.} {\it Let $E$ be a finite dimensional $k$--algebra endowed with
a $k$--linear involution $\sigma : E \to E$. Let $R_E$ be the radical of $E$, and set $\overline E = E/R_E$. Suppose that all the
orthogonal components of $\overline E$ are split. Then the canonical map $H(E,\sigma) \to \prod_{v \in \Sigma_k} H(E_v,\sigma_v)$
is injective.}

\bigskip
\noindent
{\bf Proof.} 
\noindent {\it The case of a simple algebra.}  
 Suppose first that $E$ is a simple $k$--algebra. Let $K$ be the center of $E$, and let $F$ be the fixed field of $\sigma$ in $K$.
Let $\Sigma_F$ denote the set of all places of $F$. For all $v \in \Sigma_k$, set
$E_v = E \otimes_k k_v$, and note that $E_v = \prod_{w | v} E_w$, therefore $\prod_{v \in \Sigma_k} H(E_v,\sigma_v) =
\prod_{w \in \Sigma_F} H(E_w,\sigma_w)$. By definition, $H(E,\sigma)$ is the
set of isomorphism classes of one dimensional hermitian forms over $E$. Moreover,
if $\sigma$ is orthogonal, then the hypothesis implies that $E$ is split, in other
words we have $E \simeq M_n(F)$.  Therefore the conditions
of [R 11, th. 3.3.1] are fulfilled, hence the Hasse principle holds for hermitian
forms over $E$ with respect to $\sigma$. This implies that the canonical map
$H(E,\sigma) \to \prod_{v \in \Sigma_k} H(E_v,\sigma_v)$
is injective.

\medskip
\noindent
{\it The case of a semi--simple algebra.} 
Suppose now that $E$ is semi--simple. Then 

\quad $$E \simeq E_1 \times \dots \times E_r \times
A \times A^{{\rm op}}, $$

\noindent where $E_1, \dots, E_r$ are simple algebras which are stable under the involution
$\sigma$, and where the restriction of $\sigma$ to $A \times A^{{\rm op}}$ exchanges the two
factors. Applying [BPS 13,  lemmas 3.7.1 and 3.7.2] we are reduced to the case where
$E$ is a simple algebra, and we already know that the result is true in this case.
\medskip
\noindent
{\it General case.}  
 We have $\overline E = E/R_E$. Then $\overline E$ is semi--simple,
and $\sigma$ induces a $k$--linear involution $\overline \sigma : \overline E \to \overline E$.
We have
the following commutative diagram
$$\matrix {H(E,\sigma) &  {\buildrel f \over \longrightarrow} & \prod_{v \in \Sigma_k} H(E_v,\sigma) \cr 
\downarrow & {} & \downarrow \cr 
H§(\overline E, \overline \sigma) 
& {\buildrel \overline f \over  \longrightarrow}  & \prod_{v \in \Sigma_k} H(\overline E_v,\overline \sigma), \cr }$$ 
where
the vertical maps are induced by the projection $E \to \overline E$. By [BPS 13, lemma 3.7.3], these
maps are bijective. As $\overline E$ is semi--simple,  the map $\overline f$ is injective, 
hence $f$ is  also injective. This concludes the proof. 
\medskip
\noindent
{\bf Proof of  th. 2.1.} It is clear that if $q \simeq_G q'$ over $k$, then $q \simeq_G q'$
 over all the completions of $k$. Let us prove the converse. Let $(V,h)$ be the
 $k[G]$-hermitian space corresponding to $(V,q)$, and let $\sigma : E \to E$
 be the involution induced by $(V,h)$ as in \S 1. 
Let $(V,h')$ be  the $k[G]$--hermitian space corresponding to $(V,q')$, and set 
$u = h^{-1} h'$. Then $u \in E^0$, and by lemma 1.1. the element
$[u] \in H(E,\sigma)$ determines the isomorphism class of $(V,q')$; in other
words, we have $q \simeq_G q'$ if and only if $[u] = [1]$ in $H(E,\sigma)$. Hence
the theorem is a consequence of proposition 2.2.

\bigskip
\noindent
{\bf Corollary 2.3} {\it Suppose that ${\rm char}(k)  = p > 0$, and let $V$ be a projective
$k[G]$--module.. Let
 $(V,q)$, $(V,q')$  be two $G$--forms. Then $q \simeq_G q'$ over $k$ if and only if $q \simeq_G q'$
 over all the completions of $k$.}
 
 \medskip
 \noindent
 {\bf Proof.} Since $V$ is projective, there exists a $k[G]$--module $W$ and $n \in {\bf N}$ 
 such that $V \oplus W \simeq k[G]^n$. The endomorphism ring of $k[G]^n$ is
 $M_n(k[G])$, and as ${\rm char}(k) = p  > 0$, we have $k[G] = {\bf F}_p[G] \otimes_{{\bf F}_p}k$.
 Hence $M_n(k[G])$ is isomorphic to 
 $M_n({\bf F}_p[G]) \otimes_{{\bf F}_p} k$.  Let $E = {\rm End}(V)$, let  $R_E$ be the radical of $E$, and let $\overline E = E/R_E$.
 Let us show that all the components of $\overline E$ are split. Let $e$ be the idempotent
 endomorphism of $V \oplus W$ which is the identity of $V$. 
 Set $\Lambda = {\rm End}(V \oplus W)$ and let $R_{\Lambda}$ be the radical of $\Lambda$. Then $e \Lambda e
 = E$ and $e R_{\Lambda} e = R_E$.  Set $\overline \Lambda = \Lambda / R_{\Lambda}$, and
 and let $\overline e$ be the image of $e$ in $\overline \Lambda$. Set  $\overline {k[G]} = k[G]/{\rm rad}(k[G])$. Then
 we have  $\overline E \simeq \overline e \overline \Lambda \overline e
 \simeq \overline e M_n(\overline {k[G]}) \overline e$. This implies that $\overline E$ is a component
 of the semi--simple algebra $M_n(\overline {k[G]})$. 
 Let us show that all the components of  $ M_n(\overline {k[G]})$ are split.  
 As ${\bf F}_p$ is a finite field, ${\bf F}_p[G]/({\rm rad}({\bf F}_p[G])$ is a product of matrix algebras over finite fields.
 Moreover, for any finite field $F$ of characteristic $p$, the tensor product $F \otimes_{{\bf F}_p}  k$ is a product of fields. This shows that  $({\bf F}_p[G]/({\rm rad}({\bf F}_p[G])) \otimes _{{\bf F}_p}  k$ is a product of
 matrix algebras over finite extensions of $k$; in particular, it is semi--simple. 
The natural isomorphism ${\bf F}_p[G] \otimes_{{\bf F}_p} k \to k[G]$  induces an isomorphism
 $[{\bf F}_p[G] / ({\rm rad} ({\bf F}_p[G] )) ]   \otimes_{{\bf F}_p} k \to k[G]/  ({\rm rad} ({\bf F}_p[G] ).k[G])$. Therefore
${\rm rad} ({\bf F}_p[G] .k[G])$ is the radical of $k[G]$, and we have an isomorphism 
 $[{\bf F}_p[G] / ({\rm rad} ({\bf F}_p[G] )) ]   \otimes_{{\bf F}_p} k \to k[G]/  ({\rm rad} (k[G]))$.
 Hence all the components of $k[G]/  ({\rm rad} (k[G]))$ are split. This implies that all
 the components of $\overline E$ are split as well. 
 Therefore the corollary follows from th. 2.1.
 
 \bigskip
 The following corollary is well--known (see for instance [R 11, 3.3.1 (b)]).
 \bigskip
 \noindent
 {\bf Corollary 2.4}  {\it Suppose that ${\rm char}(k) = p > 0$, and that the order of $G$ is
 prime to $p$. Then two $G$--quadratic forms are isomorphic over $k$ if and only if they
 become isomorphic over all the completions of $k$.}
 
 \medskip
 \noindent
 {\bf Proof.} This follows immediately from cor. 2.3.
 
 \bigskip
 \bigskip
 {\bf \S 3. Counter--examples to the Hasse principle}
 
 \bigskip
 Let $k$ be a  field of characteristic $p > 0$, let $C_p$ be the cyclic group
 of order $p$, and let $G = C_p \times C_p \times C_p$.
 In this section we give counter--examples to the Hasse principle for $G \times G$--quadratic
 forms over $k$ in the case where $k$ is a global field. We start with some
 constructions that are valid for any field of positive characteristic.
 
 \bigskip
 {\bf 3.1 A construction}
 
 \bigskip
 Let $D$ be a division algebra over $k$. It is well--known that there exist
 indecomposable $k[G]$--modules such that their endomorphism ring modulo
 the radical is isomorphic to $D$. We recall here such a constuction, brought to
 our attention by R. Guralnick, in order
 to use it in  3.2 in the case of quaternion algebras.
 
 \medskip
 The algebra $D$ can be generated by two elements (see for instance  [J 64, Chapter VII, \S 12, th. 3, p. 182]). Let
 us choose $i, j \in D$ be two such elements. Let us denote by $D^{\rm op}$ the
 opposite algebra of $D$, and let $d$ be the degree of $D$. 
 Then we have $D \otimes_k D^{\rm op} \simeq M_{d^2}(k)$.  Let us choose
 an isomorphism $f : D \otimes_k D^{\rm op} \simeq M_{d^2}(k)$, and set
 $a_1 = f(1 \otimes 1) = 1$,  $a_2 = f(i \otimes 1)$ and $a_3 = f(j \otimes1)$. 
 \medskip
 
  Let $g_1, g_2, g_3 \in G$ be three elements of order $p$ such that the set $\{g_1,g_2,g_3 \}$
 generates $G$.and let us define a representation $G \to {\rm GL}_{2d^2}(k)$ 
 by sending $g_m$ to the matrix

 $$ \left ( \matrix {I & a_m  \cr
0 & I  \cr} \right ) $$ for all $m = 1,2,3$. Note that this is well--defined because
${\rm char}(k) = p$. This endowes $k^{2d^2}$ with a structure of $k[G]$--module. Let
us denote by $N$ this $k[G]$--module, and let
$E_N$ be its endomorphism ring. Then 

$$E_N  = \left\{  \left ( \matrix {x & y  \cr
0 & x  \cr} \right ) \ | \ x  \in D^{\rm op} \subset M_{d^2}(k), \ y \in M_{d^2}(k)  \right \},$$
and its radical is 

$$R_N = \left\{  \left ( \matrix {0 & y  \cr
0 & 0  \cr} \right ) \ | \ y \in M_{d^2}(k)  \right \},$$
hence $E_N/R_N \simeq D^{\rm op}$.  
 
 \bigskip
 {\bf 3.2. The case of a quaternion algebra}
 
 \bigskip Let $H$ be a quaternion algebra over $k$. Then by 3.1, we get 
 a $k[G]$--module $N = N_H$  with endomorphism ring $E_N$ such
 that $E_N/R_N\simeq H^{\rm op}$, where $R_N$ is the radical of $E_N$. We now construct
 a $G$--quadratic form $q$ over $N$ in such a way that the involution it induces 
 on  $E_N/R_N \simeq H^{\rm op}$ is the canonical involution.
 
 \bigskip
Let $i, j \in H$ such that $i^2, j^2  \in k^{\times}$ and that $ij = - ji$.  Let
$\tau : H \to H$ be the orthogonal involution of $H$ obtained by
composing the canonical involution of $H$ with ${\rm Int}(ij)$. Let 
$\sigma : H^{\rm op} \to H^{\rm op}$ be the canonical involution
of $H^{\rm op}$. Let us consider the tensor product of algebras with
involution 

$$(H,\tau) \otimes (H^{\rm op},\sigma) = (M_4(k),\rho).$$ Then $\rho$ is
a symplectic involution of $M_4(k)$ satisfying  $\rho(a_m) = a_m$ for all $m = 1,2,3$,
since  $\tau(i) = (ij)(-i)(ij)^{-1} = i$,  $\tau(j) = (ij)(-j)(ij)^{-1} = j$. Let $\alpha \in M_4(k)$ be
a skew--symmetric matrix such
that for all $x \in M_4(k)$, we have $\rho(x) = \alpha^{-1} x^T \alpha$, where $x^T$
denotes the transpose of $x$. Set $A = \left ( \matrix {0 & \alpha \cr -\alpha & 0 \cr} \right )$.
Then $A^T = A$. Let $q : N \times N \to k$ be the symmetric bilinear form defined by $A$ :
$$q(v,w) = v^T A w$$ for all $v,w \in N$.  Let
$\gamma : M_8(k) \to M_8(k)$ be the involution adjoint to $q$, that is

$$\gamma(X) = A^{-1} X^T A$$ for all $X \in M_8(k)$, i.e. $q(fv,w) = q(v, \gamma(f)w)$ for
all $f \in M_8(k)$ and all $v,w \in N$. The involution $\gamma$ restricts to an involution of  $E_N$, as for all $x,y \in M_4(k)$, we have
$$\gamma  \left ( \matrix {x & y \cr 0 & x \cr} \right )  =
\left ( \matrix { \alpha^{-1}  x^T \alpha   & - \alpha^{-1}  y^T \alpha \cr 0 & \alpha^{-1}  x^T \alpha \cr }
\right ) .$$ It also sends $R_N$ to itself, and induces an involution
$\overline \gamma$ on $H^{\rm op} \simeq E_N/R_N$ that coincides with the canonical
involution of $H^{\rm op}$.

\bigskip
We claim that $q : N \times N  \to k$ is a $G$--quadratic form. To check this, it suffices to show
that $q(g_mv,g_mw) = q(v,w)$ for all $v,w \in N$ and for all $m= 1,2,3$.  Since $\rho (a_m) = a_m$
for all $m = 1,2,3$, we have

$$ \gamma  \left ( \matrix {I & a_m  \cr
0 & I  \cr} \right ) =  \left ( \matrix {I & a_m  \cr
0 & I  \cr} \right )^{-1}$$ 
and hence 
$$q(g_mv,g_mw) = q(v,\gamma(g_m)g_mw) = q(v,w)$$for all $m = 1,2,3$ and all $v,w \in N$. Thus $q$ is
a $G$--quadratic form, and by construction, the involution of $E_N$ induced by $q$ is 
the restriction of $\gamma$ to $E_N$. 

\bigskip
{\bf 3.3. Two quaternion algebras}

\medskip
Let $H_1$ and $H_2$ be two quaternion algebras over $k$. By the construction
of 3.2, we obtain two indecomposable $k[G]$--modules $N_1$ and $N_2$. Set
$E_1 = E_{N_1}$ and $E_2 = E_{N_2}$. Let $R_i$ be the radical of $E_i$ for
$i = 1,2$, and set $\overline E_i = E_i/R_i$. We also obtain $G$--quadratic
spaces  $q_i : N_i \times N_i \to k$ inducing involutions $\gamma_i : E_i \to E_i$
such that the involutions $\overline \gamma_i : \overline E_i \to \overline E_i$
coincide with the canonical involution of $H_i^{\rm op}$, for all $i = 1,2$. 
\medskip
Let us consider the tensor product $(N,q) = (N_1,q_1) \otimes_k (N_2,q_2)$. Then
$(N,q)$ is a $G \times G$--quadratic space. Set $E = {\rm End}_{k[G \times G]}(N_1 \otimes N_2)$.
Then $E \simeq E_1 \otimes E_2$. Let $I$ be the ideal of $E$ generated by
$R_1$ and $R_2$. Then there is a natural isomorphism $f : E_1 \otimes E_2 \to E$ with
$f(I) = R_E$, where
$R_E$ is the radical of $E$. Set $\overline E =
E/R_E$. Then $ \overline E  \simeq \overline E_1 \otimes \overline E_2
 \simeq H_1^{\rm op} \otimes H_2^{\rm op}$.

\bigskip
Set $\gamma = \gamma_1 \otimes \gamma_2$. Then $\gamma : E \to E$ is the
involution induced by the $G \times G$--quadratic space $(N,q)$. We obtain an
involution  $\overline \gamma
: \overline E \to \overline E$, and $\overline \gamma = \overline \gamma_1 \otimes 
\overline \gamma_2$. Let us recall that $\overline E_i = H_i^{\rm op}$ for $i = 1,2$, and
that $\overline \gamma_i$ is the canonical involution of $H_i^{\rm op}$. Hence
$\overline \gamma : \overline E \to \overline E$ is an orthogonal involution.

\bigskip
{\bf 3.4. A counter--example to the Hasse principle}
\bigskip
Suppose now that $k$ is a global field of characteristic $p$, with $p > 2$, and suppose
that $H_i$ is ramified at exactly two places $v_i, v_i'$ of $k$, such that $v_1, v_1', v_2, v_2'$
are all distinct. We have $H_1^{\rm op} \otimes H_2^{\rm op} \simeq M_2(Q)$ where
$Q$ is a quaternion division algebra over $k$, and $Q$ is ramified exactly at the places
$v_1, v_1', v_2, v_2'$ of $k$. Recall that the involution $\overline \gamma : M_2(Q) \to
M_2(Q)$ is the tensor product of the canonical involutions of $H_i^{\rm op}$. In
particular, $\overline \gamma$ is of orthogonal type.
Note that at all $v \in \Sigma_k$, one of the algebras  $H_1^{\rm op}$ or $H_2^{\rm op}$ is split.
This implies that at all $v \in \Sigma_k$, the involution $\overline \gamma$ is hyperbolic.
\medskip Let $\delta : Q \to Q$ be an orthogonal involution of the division algebra $Q$. 
Then $\overline \gamma$ is induced by some hermitian space $h : Q^2 \times Q^2 \to Q$
with respect to the involution $\delta$. As for all $v \in \Sigma_k$, the involution  $\overline \gamma$ is hyperbolic at $v$, 
the hermitian form $h$ is also hyperbolic at $v$. 
By lemma 1.1  the set of isomorphism classes of
hermitian spaces on $Q^2$ is in bijection with the set $H(\overline E,\overline \gamma)$,
the hermitian space $(Q^2,h)$ corresponding to the element $[1] \in H(\overline E,\overline \gamma)$. 

\bigskip
Let $(Q^2,h')$ be a hermitian space which becomes isomorphic to $(Q^2,h)$ over
$Q_v$ for all $v \in \Sigma_k$, but is not isomorphic to $(Q^2,h)$ over $Q$ (this is
possible by [Sch 85,  10.4.6]).
Let $u \in \overline E^0$ such that $[u] \in H(\overline E,\overline \gamma)$ corresponds
to $(Q^2,h')$ by the bijection of lemma 1.1. Then $[u] \not = [1] \in H(\overline E,\overline \gamma)$,
and the images of $[u]$ and $[1]$ coincide in $\prod_{v \in \Sigma_k}H(\overline E_v,\overline 
\gamma)$. 

\bigskip
Recall that $H(E,\gamma)$ is in bijection with the isomorphism classes of
$(G \times G)$--quadratic forms over $N$, the element $[1] \in H(E,\gamma)$ corresponding
to the isomorphism class of $(N,q)$. Let  $\pi : E \to \overline E$ be the projection,
and let  $\tilde u \in E^0$ be such that $\pi (\tilde u) = u$ (cf. lemma 1.1). Let
$(N,q')$ be a $(G \times G)$--quadratic form corresponding to $\tilde u$. The diagram
$$\matrix {H(E,\gamma) &  {\buildrel f \over \longrightarrow} & \prod_{v \in \Sigma_k} H(E_v,\gamma) \cr 
\downarrow & {} & \downarrow \cr 
H§(\overline E, \overline \gamma) 
& {\buildrel \overline f \over  \longrightarrow}  & \prod_{v \in \Sigma_k} H(\overline E_v,\overline \gamma), \cr }$$ 
is commutative, and the vertical maps are bijective by  [BPS 13, lemma 3.7.3]. Hence
$(N,q)$ and $(N,q')$ are become isomorphic over all the completions of $k$, but are
not isomorphic over $k$. 

\bigskip
\bigskip
\noindent
{\bf Bibliography}
\medskip
\noindent
[BP 11] E. Bayer--Fluckiger, R. Parimala, Galois algebras, Hasse principle and
induction--restriction methods, {\it Documenta Math.} {\bf 16} (2011), 677--707.
\medskip
\noindent
[BP 13] E. Bayer--Fluckiger, R. Parimala, Correction to : Galois algebras, Hasse principle and
induction--restriction methods, {\it Documenta Math.} {\bf 16} (2011), 677--707 (to appear).
\medskip
\noindent
[BPS 13] E. Bayer--Fluckiger, R. Parimala and J-P. Serre, Hasse principle for $G$--trace
forms, {\it Izvestia Math.}  (2013), to appear.
\medskip
\noindent
[J 64] N. Jacobson, {\it Structure of rings}, AMS Colloquium Series {\bf XXXVII} (1964).
\medskip
\noindent
[KMRT 98] M. Knus, A. Merkurjev, M. Rost \& J--P. Tignol, {\it The Book of Involutions}, AMS Colloquium Publications {\bf 44}, 1998. 
\medskip
\noindent
[M 86] J. Morales, Integral bilinear forms with a group action, {\it J. Algebra} {\bf 98} (1986), 470 - 484. 
\medskip
\noindent
[R 11] C.R. Riehm, {\it Introduction to Orthogonal, Symplectic and Unitary Representations of Finite Groups}, Fields Institute Monographies 28, AMS, 2011.
\medskip \noindent
[Sch 85] W. Scharlau, {\it Quadratic and Hermitian Forms}, Grundlehren der Math. Wiss,
Springer--Verlag (1985).

\bigskip
\bigskip

Eva  Bayer--Fluckiger

\'Ecole Polytechnique F\'ed\'erale de Lausanne

EPFL/FSB/MATHGEOM/CSAG

Station 8

1015 Lausanne, Switzerland

\medskip
eva.bayer@epfl.ch

\bigskip
\bigskip
Nivedita Bhaskhar and Raman Parimala

Department of Mathematics $ \&$ Computer Science

Emory University

Atlanta, GA 30322, USA.

\medskip
nbhaskh@emory.edu

parimala@mathcs.emory.edu

\bye